\documentclass[a4paper,11pt]{article}
\usepackage{multirow}
\usepackage{color}
\usepackage{latexsym}
\usepackage{amssymb}
\usepackage{amsfonts}
\usepackage{amsmath}
\usepackage{amsthm}
\usepackage{indentfirst}
\usepackage{mathrsfs}
\usepackage{verbatim}

\newtheorem {theorem}{Theorem}[section]

\newtheorem {lemma}{Lemma}[section]
\newtheorem {proposition}{Proposition}[section]
\newtheorem {corollary}{Corollary}[section]

\theoremstyle{remark}

\usepackage[latin1]{inputenc}

\usepackage{fullpage}

\author{Luca Ferrari\thanks{Dipartimento di Matematica e Informatica ``U. Dini", University of Firenze, Firenze, Italy\quad \texttt{luca.ferrari@unifi.it}. Partially supported by INdAM -- GNCS project CUP\_E53C23001670001.}\and Francesco Verciani\thanks{Institut f\"ur Mathematik, Universit\"at Kassel, 34132 Kassel, Germany, \, \texttt{francesco.verciani@uni-kassel.de}}}

\title{On the enumeration of permutation-invariant and complete	Naples parking functions}

\date{}

\begin{document}
	
\maketitle
	
\begin{abstract}
Naples parking functions were introduced as a generalization of classical parking functions, in which cars are allowed to park backwards, by checking up to a fixed number of previous slots, before proceedings forward as usual.
In our previous work (2024) we have provided a characterization of Naples parking functions in terms of the new notion of \emph{complete parking preference}. Our result also allowed us to describe a new characterization of permutation-invariant Naples parking functions, equivalent (but much simpler) to the one given by Carvalho et al.(2021) but using a completely different approach (and language).

In the present article we address some natural enumerative issues concerning the above mentioned objects. We propose an effective approach to enumerate permutation-invariant Naples parking functions and complete Naples parking functions which is based on some natural combinatorial decompositions. We thus obtain formulas depending on some (generally simpler) quantities, which are of interest in their own right, and that can be described in a recursive fashion.  
\end{abstract}
	
\section{Introduction}

The notion of parking function was introduced by Konheim and Weiss \cite{KW} in a computer science context, specifically in connection with a hashing problem. Since then, parking functions showed up in several other mathematical contexts, such 
as hyperplane arrangements \cite{St}, representation theory \cite{KT} and combinatorial geometry \cite{AW}. Roughly speaking, a parking function describes the preferences of a set of cars that successfully park in a series of parking slots arranged in linear order in a one-way street. Each car has a preferred slot and, if this is already occupied, the car is allowed to park into the next available slot. A more formal and precise definition of parking function is given in Section \ref{prel}.

In recent years several generalizations of parking functions have been proposed, each of which describes a different parking model \cite{CCHJR}. Among them, the most investigated one is probably that of Naples parking function, in which cars are required to check a limited number of parking slots before its preferred one (in case this is already occupied).

The combinatorics of Naples parking functions appear to be significantly more complicated than the combinatorics of classical parking functions. For instance, a rearrangement of a Naples parking function may not be a Naples parking function, in contrast to the classical case. As a consequence, while most of the combinatorics of classical parking functions is rather well understood, the same cannot be said for Naples parking functions. In particular, there are few enumerative results regarding Naples parking functions, and most of them concern specific classes. Some interesting recursive formulas for the total number of $k$-Naples parking of length $n$ have been found, however a  satisfactory closed formula is still lacking; the only general result in this direction is contained in \cite{CHJKRSV}, where a very complicated expression in terms of a sum over permutations is presented.

In the present work we build on what we did in our previous paper \cite{FV} to obtain some new enumerative results on Naples parking functions. Using the tools developed in \cite{FV}, and recalled in Section \ref{prel} together with all standard definitions and notations regarding classical and Naples parking functions, we are able to effectively enumerate the class of permutation-invariant Naples parking functions (Section \ref{enum_perm_inv}) and the class of complete Naples parking functions (Section \ref{enum_compl}), in both cases with respect to the length. More precisely, we describe an effective counting procedure to enumerate permutation-invariant Naples parking functions, thanks to which we are able to produce several of such counting coefficients, whereas in the case of complete Naples parking functions we determine a recursion which, again, allows us to compute the desired values very efficiently.  Most of the coefficients we find are not recorded in \cite{Sl}, with a few intriguing exceptions. For a couple of those we are able to prove the conjectured hit.

\section{Preliminary definitions and results}\label{prel}

A \emph{parking preference} $\alpha =(a_1 ,a_2 ,\ldots ,a_n)$ of length $n$ is a $n$-tuple of positive integers between 1 and $n$ to be interpreted as follows: for an ordered list of $n$ cars attempting to park on $n$ numbered slots in a street, $a_i$ denotes the preferred slot for the $i$-th car. The set of all parking preferences of length $n$ will be denoted $PP_n$. The total number of cars having preference $i$ will be denoted with $|\alpha |_i$. The \emph{classical} parking rule requires that each car tries to park into its preferred slot and, if this is already occupied, then proceeds forward and parks into the first available slot, if any, otherwise it simply exits the street (thus failing to park). If all cars successfully park using the above rule, then we say that $\alpha$ is a \emph{parking function} of length $n$, and we denote the set of all those with $PF_n$. An elegant characterization of parking functions says that $\alpha$ is a parking function if and only if $\sum_{i=1}^{j}|\alpha |_i \geq j$, for all $j\leq n$. The latter condition is equivalent to the fact that the nondecreasing rearrangement $(b_1 ,b_2 ,\ldots ,b_n)$ of $\alpha$ is such that $b_j \leq j$, for all $j\leq n$. As a consequence, we get that any rearrangement of a parking function is a parking function as well.

The classical parking rule can be modified into the \emph{$k$-Naples parking rule}, by allowing each car to check up to $k$ slots before its preferred one, and then proceed forward if it does not find an available place. A \emph{$k$-Naples parking function} is a parking preference such that all cars succeed to park following the $k$-Naples parking rule, and the set of $k$-Naples parking functions of length $n$ will be denoted $PF_{n,k}$. We remark that, unlike classical parking function, a rearrangement of a Naples parking function is not a Naples parking function in general.

\bigskip

We now recall the main notations and facts we need from \cite{FV}.

Let $\alpha =(a_1 ,a_2 ,\ldots,a_n )\in PP_n$. Then, for every $j\leq n$, we set
\[
u_\alpha (j)=\sum_{i=j}^{n}|\alpha|_i -(n-j+1).
\]

In other words, $u_\alpha (j)\in \mathbb{Z}$ denotes the excess of cars having preference at least $j$ with respect to the number of available slots from $j$ onwards. The map $u_{\alpha}$ is called the \emph{excess function} of the parking preference $\alpha$. We set $U_{\alpha}=\{ j\leq n\, |\, u_\alpha (j)\geq 1\}$ and we remark that obviously $1\notin U_{\alpha}$, for all parking preferences $\alpha$. The set $U_{\alpha}$ is called the \emph{excess set} of $\alpha$. A parking preference $\alpha$ (of length $n$) such that $U_{\alpha}=[2,n]$ (that is, $U_{\alpha}$ is maximum) is said to be \emph{complete}.

In the next couple of lemmas we record some basic properties of the excess function.

\begin{lemma}\label{elementary_properties}
	Let $\alpha \in PP_n$. For all $j,j_1 ,j_2 \leq n$, with $j_1 <j_2$, we have: %(i) $u_\alpha (j)=j-1-\sum_{i=1}^{j-1}|\alpha |_i$; (ii) $u_\alpha (j_2 )-u_\alpha (j_1 )=j_2 -j_1 -\sum_{i=j_1 }^{j_2 -1}|\alpha |_i$; (iii) $u_\alpha (j)=u_\alpha (j+1)+|\alpha |_j -1$, when $j<n$; (iv) $u_\alpha (j)<j$.
	\begin{itemize}
		\item[(i)] $u_\alpha (j)=j-1-\sum_{i=1}^{j-1}|\alpha |_i$;
		\item[(ii)] $u_\alpha (j_2 )-u_\alpha (j_1 )=j_2 -j_1 -\sum_{i=j_1 }^{j_2 -1}|\alpha |_i$;
		\item[(iii)] $u_\alpha (j)=u_\alpha (j+1)+|\alpha |_j -1$, when $j<n$;
		%\item[(iv)] $u_\alpha (j)<j$.
	\end{itemize} 
\end{lemma}

\begin{lemma}\label{elementary_intervals}
	Let $\alpha \in PP_n$ and let $[p,q]$ be a maximal interval of consecutive positions contained in $U_\alpha$. Then we have: %(i) $u_\alpha (p)=1$; (ii) $u_\alpha (p-1)=0$; (iii) $|\alpha |_{p-1}=0$; (iv) $|\alpha |_q \geq 2$. 
	\begin{itemize}
		\item[(i)] $u_\alpha (p)=1$;
		\item[(ii)] $u_\alpha (p-1)=0$;
		\item[(iii)] $|\alpha |_{p-1}=0$;
		\item[(iv)] $|\alpha |_q \geq 2$.
	\end{itemize}     
\end{lemma}

Using Lemma \ref{elementary_properties}, it is not difficult to show that $\alpha \in PF_n$ if and only if $U_\alpha =\emptyset$. Such a result can be generalized to a necessary condition for $k$-Naples parking function, which is reported below in terms of the excess function.

\begin{proposition}\label{necessary_PF}
	If $\alpha \in PF_{n,k}$, then $u_\alpha (j)\leq k$, for all $j\leq n$.
\end{proposition}

Given $\alpha =(a_1 ,a_2 ,\ldots ,a_n ) \in PP_n$ and $w\in \mathbb{Z}$, the \emph{$w$-shift} of $\alpha$ is the $n$-tuple $\tau_w (\alpha )=(a_1 -w, a_2 -w,\ldots ,a_n -w)$. We also extend the above notation to sets: given $W\subseteq \mathbb{Z}$, we set $\tau_w (W)=\{ i\in \mathbb{Z}\, |\, i+w\in W\}$. Moreover, if $J=\{ j_1 ,j_2 ,\ldots j_h \} \subseteq [1,n]$, with $j_1 <j_2 <\cdots <j_h$, the \emph{restriction} of $\alpha$ to $J$ is the $h$-tuple $\alpha_{|_{J}}=(a_{j_1},a_{j_2},\ldots ,a_{j_h})$. The following lemma collects some useful facts concerning shifts and restrictions of parking preferences.

\begin{lemma}\label{shifts_and_restrictions}
	Let $\alpha \in PP_n$ and let $j\leq n$ such that $u_{\alpha}(j)=0$. Setting $J=\{ i\leq n\, |\, a_i \geq j\}$ (and $J^c =[1,n]\setminus J$), we have that $|J|=n-j+1$, hence $\alpha_{|_{J^c}}\in PP_{j-1}$ and $\tau_{j-1}(\alpha_{|_{J}})\in PP_{n-j+1}$ (and, if $\alpha$ is a $k$-Naples parking function, then $\tau_{j-1}(\alpha_{|_{J}})$ is a $k$-Naples parking function as well). Moreover $u_{\alpha_{|_{J^c}}}(i)=u_{\alpha}(i)$ for all $i\leq j-1$, hence $U_{\alpha_{|_{J^c}}}=U_{\alpha}\cap [1,j-1]$, and similarly $u_{\tau_{j-1}(\alpha_{|_{J}})}(i)=u_{\alpha}(i+j-1)$ for all $i\leq n-j+1$, hence $U_{\tau_{j-1}(\alpha_{|_{J}})}=\tau_{j-1}(U_{\alpha}\cap [j,n])$. 
\end{lemma}

The following characterization of permutation-invariant $k$-Naples parking functions is one of the main results of \cite{FV}.

\begin{proposition}\label{per_inv_k}
	Let $\alpha =(a_1 ,a_2 ,\ldots ,a_n ) \in PP_n$ and $k\geq 1$. The following are equivalent:
	\begin{itemize}
		\item[(i)] for every maximal interval $[p,q]\subseteq U_{\alpha}$, $|[p,q]|\leq k$;
		\item[(ii)] for all permutations $\sigma$ of size $n$, the $n$-tuple $\sigma (\alpha )=(a_{\sigma (1)} ,a_{\sigma (2)} ,\ldots ,a_{\sigma (n)})$ is a k-Naples parking function of length $n$. In particular, $\alpha \in PF_{n,k}$.
	\end{itemize}
\end{proposition}

\section{Permutation-invariant Naples parking functions}\label{enum_perm_inv}

Denote with $T_{n,k}$ the number of permutation-invariant $k$-Naples parking functions of length $n$ and with $t_{n,k}$ the number of those for which $2\in U_{\alpha}$. Our first goal is to express $T_{n,k}$ in terms of the simpler coefficients $t_{n,k}$.

\begin{proposition}\label{T_and_t}
	For all $n,k\geq 0$, and setting $t_{0,k}=1$, we have
	\begin{equation}\label{T_vs_t}
	T_{n,k}=\sum_{j=0}^{n}{n\choose j}(j+1)^{j-1}t_{n-j,k}.
	\end{equation}
\end{proposition}

\begin{proof}
	We start by observing that, using Proposition \ref{necessary_PF}, it is easy to realize that every parking function is a $k$-Naples parking function (for any $k$), and so in particular permutation-invariant. Therefore we will focus on Naples parking functions that are \emph{not} parking functions.
	So let $\alpha =(a_1 ,a_2 ,\ldots ,a_n) \in PF_{n,k}\setminus PF_n$ be permutation-invariant, then $U_{\alpha}\neq \emptyset$, and we set $h=\min U_{\alpha}\geq 2$ and $H=\{ i\leq n\, |\, a_i \geq h-1\}$. Lemma \ref{elementary_intervals} implies that $|H|\geq 2$ and that $u_{\alpha}(h-1)=0$.
	
	Thanks to Lemma \ref{shifts_and_restrictions}, we have that $|H^c|=h-2$, $\alpha_{|_{H^c}}\in PP_{h-2}$ and $U_{\alpha_{|_{H^c}}}=U_{\alpha}\cap [1,h-2]=\emptyset$, hence $\alpha_{|_{H^c}}\in PF_{h-2}$. Such a lemma also implies that $\tau_{h-2}(\alpha_{|_{H}})\in PP_{n-h+2}$ and $U_{\tau_{h-2}(\alpha_{|_{H}})}=\tau_{h-2}(U_{\alpha}\cap [h-1,n])=\tau_{h-2}(U_{\alpha})$. Since $\alpha$ is $k$-Naples permutation-invariant, all maximal intervals of $U_{\alpha}$ have size at most $k$ (by Proposition \ref{per_inv_k}). The same holds for $U_{\tau_{h-2}(\alpha_{|_{H}})}$ and so $\tau_{h-2}(\alpha_{|_{H}})\in PF_{n-h+2,k}$ and is permutation-invariant. Finally, we also have that $u_{\tau_{h-2}(\alpha_{|_{H}})}(2)=u_{\alpha}(h)=1$ (by Lemma \ref{elementary_intervals}), thus $2\in U_{\tau_{h-2}(\alpha_{|_{H}})}$.
	
	Summing up, the above analysis shows that a permutation-invariant $k$-Naples parking function of length $n$ which is not a parking function is uniquely determined by a set of positive integers $H\subseteq [n]$, with $|H|=n-h+2$, together with a parking function $\alpha_1 \in PF_{h-2}$ (on the cars associated with $H^c$) and a permutation-invariant $k$-Naples parking function $\alpha_2 \in PF_{n-h+2,k}$ such that $2\in U_{\alpha_2}$. Since $2\leq |H|\leq n$, we have that $2\leq h\leq n$. Setting $j=h-2$, we thus get that the number of permutation-invariant $k$-Naples parking functions of length $n$ which are not parking functions is given by $\sum_{j=0}^{n-2}{n\choose j}(j+1)^{j-1}t_{n-j,k}$. To conclude, we observe that, if we complete the previous summation with the summands corresponding to $j=n-1$ and $j=n$, then the former one vanishes (since clearly $t_{1,k}=0$), whereas the latter one counts (standard) parking functions of length $n$: this gives precisely formula (\ref{T_vs_t}).       
\end{proof}

The second step of our approach consists of expressing the coefficients $t_{n,k}$ in terms of somehow simpler ones. Given $n\geq 0, k\geq 1$,  we set $\Theta_{n,k}^{(=)}=|\{ \alpha \in PP_n \, |\, U_{\alpha}=[2,k+1] \}|$ and $\Theta_{n,k}^{(\leq)}=\sum_{i=1}^{k}\Theta_{n,i}^{(=)}$. In other words, $\Theta_{n,k}^{(\leq)}$ is the number of parking preferences whose excess set is a single interval of size at most $k$ having minimum 2. We observe that clearly $\Theta_{0,k}^{(=)}=\Theta_{1,k}^{(=)}=0$ for all $k\geq 1$, and also that $\Theta_{n,k}^{(=)}=0$ when $n\leq k$.

Notice that, as a consequence of Proposition \ref{per_inv_k}, all parking preferences counted by $\Theta_{n,k}^{(\leq)}$ are $k$-Naples parking functions (and the same is true for $\Theta_{n,i}^{(=)}$ when $i\leq k$). An argument analogous to that of the previous proposition provides the following recursive description of $t_{n,k}$.

\begin{proposition}\label{t_and_theta}
	For all $n\geq 0, k\geq 1$, we have
	\begin{equation}\label{t_vs_theta}
	t_{n,k}=\sum_{i=0}^{n}{n\choose i}\Theta_{i,k}^{(\leq)}t_{n-i,k}.
	\end{equation}	
\end{proposition}

\begin{proof}
	Let $\alpha \in PF_{n,k}$ be permutation-invariant, with $2\in U_{\alpha}$. Suppose that $[2,j]$ is maximal in $U_{\alpha}$, with $j\leq k+1$. Denote with $h$ the minimum of $U_{\alpha}\setminus [2,j]$ (if it exists), and set $H=\{ i\leq n\, |\, h-1\leq a_i \leq n\}$. It is clear that $\alpha_{|_{H^c}}$ is a $k$-Naples parking function such that $U_{\alpha_{|_{H^c}}}$ consists of the single interval $[2,j]$, whereas $\tau_{h-2}(\alpha_{|_{H}})$ is a permutation-invariant $k$-Naples parking function such that $2\in U_{\tau_{h-2}(\alpha_{|_{H}})}$. This gives exactly formula (\ref{t_vs_theta}).
\end{proof}

We remark that, since $\Theta_{0,k}^{(\leq)}=0$, in the right hand side of (\ref{t_vs_theta}) the summand corresponding to $i=0$ vanishes, thus giving a properly defined recursion to determine $t_{n,k}$.

Propositions \ref{T_and_t} and \ref{t_and_theta} give a recipe to express the number of permutation-invariant Naples parking functions as a function of the $\Theta_{n,k}^{(\leq)}$'s only. From a purely combinatorial point of view, this can be expressed by saying that a permutation-invariant $k$-Naples parking function $\alpha$ of length $n$ can be uniquely determined as follows. Given a composition $\beta=(b_0 ,b_1 ,\ldots ,b_m)$ of $n$ (i.e., $b_0 +b_1 +\cdots +b_m =n$), with $b_0 \geq 0$ and $b_i \geq 2$ for all $i\neq 0$, choose a set partition $\{ B_0 ,B_1 ,\ldots ,B_m\}$ of $\{ 1,2,\ldots n\}$ such that $|B_i |=b_i$, for all $i$. Then assign: 
\begin{itemize}
	\item a parking function of size $b_0$ on $B_0$ (these will be those cars having preferences between 1 and $b_0$), and
	\item for each $i=1,2,\ldots ,m$, a permutation-invariant $k$-Naples parking function $\alpha_i$ of size $b_i$ with the property that, for all $i$, the excess set $U_{\alpha_i}$ consists of a single interval (of cardinality at most $k$) containing 2 (these will be those cars having preferences between $\sum_{j=0}^{i-1}b_j +1$ and $\sum_{j=0}^{i}b_j$).   
\end{itemize} 

%Given a composition $\beta=(b_0 ,b_1 ,\ldots ,b_m)$ of $n$ (i.e., $b_0 +b_1 +\cdots +b_m =n$), with $b_i \geq 2$ for all $i\neq 0$, choose the set $U_\alpha =I_1 \cup I_2 \cup \cdots \cup I_m$, where the $I_j$'s are disjoint intervals such that $2\leq |I_j |\leq b_j$ and $|I_j |\leq k$, for all $j=1,2,\ldots ,m$. Then assign:
%\begin{itemize}
%	\item a parking function of size $b_0$ on the set of cars whose preferences are at most $\min I_1 -2$, and
%	\item a permutation-invariant $k$-Naples parking function $\alpha_j$ of size $b_j$ on each of the sets of cars having preferences between $\min I_j -1$ and $\min I_{j+1}-2$, with the property that, for all $j$, $U_{\alpha_j}$ consists of a single interval (of cardinality $|I_j |\leq k$) containing 2.   
%\end{itemize}

%the set $U_\alpha =[b_0 +1,y_1 ]\cup [b_0 +b_1 +1,y_2 ]\cup \cdots \cup [b_0 +b_1 +\cdots +b_{m-1}+1,y_m ]$ (where each interval has cardinality at most $k$) and then assigning: 
%\begin{itemize}
%	\item a parking function on the set of cars whose preferences lie between 1 and $b_0$, and
%	\item a permutation-invariant $k$-Naples parking function $\alpha_i$ on each of the sets of cars having preferences between $\sum_{j=0}^{i-1}b_j +1$ and $\sum_{j=0}^{i}b_j$, with the property that, for all $i$, $U_{\alpha_i}$ consists of a single interval (of cardinality at most $k$) containing 2.   
%\end{itemize}

The above considerations can be translated into the following counting formula
\begin{equation}\label{general_formula}
T_{n,k}=\sum_m \sum_{\beta \in W_m}\binom{n}{b_0 ,b_1, \ldots ,b_m}(b_0 +1)^{b_0 -1}\prod_{i=1}^{m}\Theta_{b_i ,k}^{(\leq)},
\end{equation}
where $W_m=\{ (b_0, b_1 ,\ldots ,b_m )\, |\, b_0 +b_1 +\cdots +b_m =n, b_0\geq 0, b_i\geq 2, i=1,2,\ldots ,m\}$ and $\binom{n}{b_0 ,b_1, \ldots ,b_m}=\frac{n!}{b_0 !b_1 !\cdots b_m!}$ denotes the usual multinomial coefficient.

\bigskip

As it is clear from formula (\ref{general_formula}), in order to effectively compute $T_{n,k}$ we need to find information on the coefficients $\Theta_{n,k}^{(\leq)}$. In the rest of the present section we will in fact provide a recursive recipe to compute the coefficients $\Theta_{n,k}^{(=)}$, by which we can in turn compute the $\Theta_{n,k}^{(\leq)}$'s. We start by showing an interesting symmetry property.

\begin{lemma}\label{symmetry}
	For any $n\geq 2$ and $1\leq k\leq n-1$, we have $\Theta_{n,k}^{(=)}=\Theta_{n,n-k}^{(=)}$. 
\end{lemma}

\begin{proof}
	Given $\alpha =(a_1 ,a_2 ,\ldots ,a_n)\in PP_n$, with $U_{\alpha}=[2,k+1]$, define $\tilde{\alpha}=(n+2-a_1 ,n+2-a_2 ,\ldots ,n+2-a_n)$ (this will be called the \emph{complement} of $\alpha$). Since $|\alpha |_1 =0$ (by (iii) of Lemma \ref{elementary_intervals}), it follows immediately that $\tilde{\alpha}\in PP_n$. By definition of $\tilde{\alpha}$, we have that $|\tilde{\alpha}|_i =|\alpha|_{n+2-i}$, hence
	\begin{align*}
	u_{\tilde{\alpha}}(j)&=\sum_{i=j}^{n}|\tilde{\alpha}|_i -(n-j+1)=\sum_{i=j}^{n}|\alpha |_{n+2-i}-(n-j+1) \\
	&=\sum_{i=1}^{n+2-j}|\alpha |_i -(n-j+2)+1=-u_{\alpha}(n-j+3)+1
	\end{align*}
	(the last equality follows from (i) of Lemma \ref{elementary_properties}). As a consequence, we get that $u_{\tilde{\alpha}}(j)\geq 1$ if and only if $u_{\alpha}(n-j+3)\leq 0$, i.e. $j\in U_{\tilde{\alpha}}$ if and only if $n-j+3\notin U_{\alpha}=[2,k+1]$, or, equivalently, if and only if $j\leq n-k+1$. Therefore we have that $U_{\tilde{\alpha}}=[2,n-k+1]$, that is $\tilde{\alpha}\in \Theta_{n,n-k}^{(=)}$. Since the complement map is clearly an involution, we get that $\Theta_{n,k}^{(=)}=\Theta_{n,n-k}^{(=)}$, as desired.
\end{proof}

Next we further decompose the quantities $\Theta_{n,k}^{(=)}$ by introducing more parameters.

Let $n,m,h\geq 1$ and $k\geq 0$. Define $\vartheta_{n,k}(m,h)$ to be the number of parking preferences $\alpha \in PP_n$ such that $U_{\alpha}=[2,k+1]$, $|\alpha |_i =0$ for $i=1,\ldots h-1$, and $|\alpha |_h =m$ (in particular, if $k=0$ we have that $U_{\alpha}=\emptyset$). Roughly speaking, the coefficients $\vartheta$ refine the eunumeration of parking preferences whose excess set is a single interval by looking at the smallest preference. More specifically
\begin{equation}\label{refine}
\Theta_{n,k}^{(=)}=\sum_{m\geq 1}\sum_{h\geq 1}\vartheta_{n,k}(m,h).
\end{equation}

The expression (\ref{refine}) can be considerably simplified by observing that many of the coefficients $\vartheta_{n,k}(m,h)$ actually vanish.

\begin{lemma}\label{vanishing}
	Let $n,m,h\geq 1$ and $n>k\geq 0$.
	\begin{itemize}
		\item[(i)] $\vartheta_{n,k}(m,h)=0$ for all $h>k+1$; in particular, $\vartheta_{n,0}(m,h)=0$ for all $h>1$;
		\item[(ii)] if $k\geq 1$, then $\vartheta_{n,k}(m,1)=0$;
		\item[(iii)] if $k\geq 1$, then $\vartheta_{n,k}(m,k+1)=0$ for all $m<k+1$;
		\item[(iv)] if $k\geq 1$ and $2\leq h\leq k$, then $\vartheta_{n,k}(m,h)=0$ for all $m\geq h$;
		\item[(v)] in all the remaining cases, $\vartheta_{n,k}(m,h)\neq 0$. 
	\end{itemize}
\end{lemma}

\begin{proof}
	Let $\alpha \in PP_n$. If $\alpha \in PF_n$, then $|\alpha |_1 \neq 0$, hence $\vartheta_{n,0}(m,h)=0$ for all $h>1$. For the rest of the proof we assume that $k\geq 1$ and that $U_{\alpha}=[2,k+1]$. By Lemma \ref{elementary_intervals}, we have that $|\alpha |_{k+1}\neq 0$, and so $\vartheta_{n,k}(m,h)\neq 0$ for all $h>k+1$, thereby proving $(i)$. Morevoer, the same lemma also gives that $|\alpha |_1 =0$, from which $(ii)$ follows. Now suppose that $h=k+1$. This means that $|\alpha |_i =0$ for all $i\leq k$, and so, using Lemma \ref{elementary_properties} and the fact that $k+2 \notin U_{\alpha}$,
	\[
	0\geq u_{\alpha}(k+2)=k+1-\sum_{i=1}^{k+1}|\alpha |_i =k+1-|\alpha |_{k+1}.
	\]
	Therefore $|\alpha |_{k+1}\geq k+1$, and this implies that $\vartheta_{n,k}(m,k+1)=0$ for all $m<k+1$, which is $(iii)$. On the other hand, if $2\leq h\leq k$, then $h+1 \in U_{\alpha}$, hence
	\[
	1\leq u_{\alpha}(h+1)=h-\sum_{i=1}^{h}|\alpha |_i =h-|\alpha |_h ,
	\] 
	i.e. $|\alpha |_h \leq h-1$, which proves $(iv)$.
	
	To conclude the proof, we now need to show that in all the remaining cases $\vartheta_{n,k}(m,h)\neq 0$. We will prove this by explicitly providing a parking preference for each single case.
	
	\begin{itemize}
		\item If $k=0$ and $h=1$, choose $\alpha =(\underbrace{1,\ldots ,1}_{m},\underbrace{2,\ldots ,2}_{n-m})$.
		\item If $k\geq 1$ and $2\leq h\leq k$, for $m<h$, choose $\alpha =(\underbrace{h,\ldots ,h}_{m},\underbrace{k+1,\ldots ,k+1}_{n-m})$.
		\item If $k\geq 1$ and $h=k+1$, for $m\geq k+1$, choose $\alpha =(\underbrace{k+1,\ldots ,k+1}_{m},\underbrace{k+2,\ldots ,k+2}_{n-m})$   
	\end{itemize}

\end{proof}

Using the above lemma, we can therefore simplify formula (\ref{refine}) by keeping only the nonzero summands, thus getting

\[
\Theta_{n,k}^{(=)}= \sum_{h=2}^{k}\sum_{m=1}^{h-1}\vartheta_{n,k}(m,h)+\sum_{m=k+1}^{n}\vartheta_{n,k}(m,k+1) .
\]

To conclude our enumeration strategy, we now give a recurrence to compute the coefficients $\vartheta$.

\begin{proposition}\label{base}
	For $n\geq m\geq 1$, we have
	\[
	\vartheta_{n,0}(m,1)=\binom{n-1}{m-1}n^{n-m}.
	\]
\end{proposition} 

\begin{proof}
	Observe that $\vartheta_{n,0}(m,1)$ counts (classical) parking functions of length $n$ in which exactly $m$ cars have preference 1. The generating function of such coefficients can be found in \cite{Y} (Corollary 1.16).  
\end{proof}

\begin{proposition}\label{step}
	Let $n\geq k\geq 1$ and $1\leq h\leq n$. Then
	\begin{equation}\label{m=1}
		\vartheta_{n,k}(1,h)=n\cdot \sum_{h'= h}^{k}\sum_{m'=1}^{n-1}\vartheta_{n-1,k-1}(m',h')
	\end{equation}
	and, for $m\geq 2$,
	\begin{equation}\label{m>=2}
		\vartheta_{n,k}(m,h)=\frac{n}{m}\vartheta_{n-1,k-1}(m-1,h-1).
	\end{equation}
\end{proposition}

\begin{proof}
	The basic idea is to develop a construction showing how each of the $k$-Naples parking functions of length $n$ counted by $\vartheta_{n,k}(m,h)$ can be obtained (not necessarily in a unique way) from one of length $n-1$.
	
	Let $\alpha \in PF_{n,k}$ such that $U_{\alpha}=[2,k+1]$, and let $h\leq n$ be such that $|\alpha |_i =0$ for all $i<h$ and $|\alpha |_h \neq0$. Choose an index $\tilde{i}$ such that $a_{\tilde{i}} =h$ and denote with $J$ the set of the remaining indices (that is, $J=\{ 1,\ldots ,n\} \setminus \{ \tilde{i}\}$). Now consider the parking preference $\tau_{1}(\alpha_{|_{J}})$ (of length $n-1$); roughly speaking, this is nothing else than removing $a_{\tilde{i}}$ from $\alpha$ and decreasing the remaining entries by 1. We now want to show that $\tau_{1}(\alpha_{|_{J}})$ is a $(k-1)$-Naples parking function of length $n-1$ counted by a suitable $\vartheta$ coefficient.
	
	Since obviously $|\tau_{1}(\alpha_{|_{J}})|_i =|\alpha|_{i+1}$ for all $i\neq h-1$, and $|\tau_{1}(\alpha_{|_{J}})|_{h-1}=|\alpha|_h -1$, an easy computation shows that, for $h\leq j\leq n$, we have $u_{\tau_{1}(\alpha_{|_{J}})}(j)=u_{\alpha}(j+1)$. This fact, together with the fact that $h\leq k+1$ (which is a consequence of (i) of Lemma \ref{vanishing}, since otherwise $\alpha$ would not exist), tells us that $[h,k]\subseteq U_{\tau_{1}(\alpha_{|_{J}})}$, and also that $i\notin U_{\tau_{1}(\alpha_{|_{J}})}$ for all $i>k$. Furthermore, an similarly easy computation as above shows that, for $2\leq j\leq h-1$, we have $u_{\tau_{1}(\alpha_{|_{J}})}(j)=u_{\alpha}(j+1)-1$. Notice that, for $3\leq i\leq h$, we have $u_{\alpha}(i)\geq 2$; in fact, if this were not the case, then the fact that $|\alpha |_{i-1}=0$ would imply (using (iii) of Lemma \ref{elementary_properties}) that $u_{\alpha}(i-1)=u_{\alpha}(i)+|\alpha |_{i-1} -1\leq 0$, against the fact that $i-1\in U_{\alpha}$ (this is because $i-1\in [2,h-1]\subset [2,k+1]=U_{\alpha}$). As a consequence, we have that, for all $2\leq j\leq h-1$, $u_{\tau_{1}(\alpha_{|_{J}})}(j)=u_{\alpha}(j+1)-1\geq 1$, and so $[2,h-1]\subseteq U_{\tau_{1}(\alpha_{|_{J}})}$. Summing up, we have thus proved that $U_{\tau_{1}(\alpha_{|_{J}})}=[2,k]$. This ensures that $\tau_{1}(\alpha_{|_{J}})$ is a $(k-1)$-Naples parking function of length $n-1$.
	
	In order to reverse the above construction, let $\beta$ be a $(k-1)$-Naples parking function of length $n-1$ such that $U_{\beta}=[2,k]$, and let $h\leq k$ be such that $|\beta |_i =0$ for all $i\leq h-2$. We can thus construct a $k$-Naples parking function $\alpha$ of length $n$ by simply increasing all entries of $\beta$ by 1 and then inserting somewhere a new entry equal to $h$. Using similar arguments as above it is possible to prove that indeed $\alpha \in PF_{n,k}$, $U_{\alpha}=[2,k+1]$, and of course $|\alpha|_i =0$ for all $i\leq h-1$. However, several different parking preferences $\alpha$ can be obtained from the same $\beta$, since there are $n$ possible positions where to insert the new entry: this explains the factor $n$ appearing in the left hand side of both (\ref{m=1}) and (\ref{m>=2}). Moreover, we need to distinguish two cases, depending on whether $\beta$ contains some entry equal to $h-1$ or not. If $|\beta |_{h-1}=0$, then the above ``reverse construction" allows us to reconstruct $\alpha$ in a unique way, so we get (\ref{m=1}). On the other hand, if $|\beta |_{h-1}=m-1$, with $m\geq 2$, then the same $\alpha$ can be obtained from $m$ different preferences $\beta$, and this gives (\ref{m>=2}).      
\end{proof}

The recursion described in Propositions \ref{base} and \ref{step} allows us to effectively compute the coefficients $\vartheta_{n,k}(m,h)$, which in turn make it possible to effectively compute the coefficients $\Theta_{n,k}^{(=)}$ (by simply using (\ref{refine}) or, even better, its simplified version displayed before Proposition \ref{base}). The table below reports the first few values of such coefficients.

\begin{center}
\begin{table}[h!]
\begin{tabular}{c|ccccccc}
$\Theta_{n,k}^{(=)}$ & $k=1$ & $k=2$ & $k=3$ & $k=4$ & $k=5$ & $k=6$ & $k=7$ \\
\hline
$n=2$ & 1      &        &        &        &        &        & \\
$n=3$ & 4      & 4      &        &        &        &        & \\
$n=4$ & 27     & 21     & 27     &        &        &        & \\
$n=5$ & 256    & 176    & 176    & 256    &        &        & \\
$n=6$ & 3125   & 1995   & 1765   & 1995   & 3125   &        & \\
$n=7$ & 46656  & 28344  & 23304  & 23304  & 28344  & 46656  & \\
$n=8$ & 823543 & 482825 & 378007 & 351337 & 378007 & 482825 & 823543 \\ 	
\end{tabular}
\caption{The first values of the coefficients $\Theta^{(=)}_{n,k}$.}\label{theta=}
\end{table}
\end{center}

As we already know from Lemma \ref{symmetry}, the rows of Table \ref{theta=} are symmetric. In general, rows and columns of such a table (as well as the whole table) seem not to be recorded in \cite{Sl}. A notable exception is the first column, i.e. the values $\Theta^{(=)}_{n,1}=\Theta^{(=)}_{n,n-1}=(n-1)^{n-1}$. In the next proposition we will prove this fact in a bijective fashion, namely by showing that the set of 1-Naples parking functions counted by $\Theta^{(=)}_{n,1}$ are in bijection with prime parking functions of length $n$. A \emph{prime parking function} is a parking function $\alpha =(a_1 ,a_2,\ldots ,a_n )$ such that the sequence obtained from $\alpha$ by removing a 1 is still a parking function. Equivalently, $\alpha$ is a prime parking function when, denoting with $\beta =(b_1 ,b_2 ,\ldots ,b_n )$ its nondecreasing rearrangement, we have that $b_i <i$ for all $i\neq 1$. According to Stanley \cite{St}, the notion of prime parking function was introduced by Gessel, who also gave the above mentioned enumerative formula (see also \cite{DO}).

\begin{proposition}
	For all $n\geq 2$, we have
	\[
	\Theta^{(=)}_{n,1}=(n-1)^{n-1}.
	\]
\end{proposition}

\begin{proof}
	Let $\alpha =(a_1 ,a_2,\ldots ,a_n )$ be a 1-Naples parking function such that $U_{\alpha}=\{ 2\}$. From Lemma \ref{elementary_intervals} we get that $|\alpha |_1 =0$ and $|\alpha |_2 \geq 2$. This implies that the sequence $\tilde{\alpha}=(a_1 -1, a_2 -1,\ldots ,a_n -1)$ is a parking preference of length $n$ which contains at least two 1's. Moreover, given $j>1$, since $j+1\notin U_{\alpha}$ (and so $u_{\alpha}(j+1)\leq 0$), we can use (i) of Lemma \ref{elementary_properties} to get that $\sum_{i=1}^{j}|\alpha|_i =j-u_{\alpha}(j+1)\geq j$. As a consequence, in the nondecreasing rearrangement of $\alpha$ the $j$-th entry is at most $j$, hence in the nondecreasing rearrangement of $\tilde{\alpha}$ the $j$-th entry is strictly less than $j$. We have thus proved that $\tilde{\alpha}$ is a prime parking function. Since the above construction is clearly reversible, it is a bijection, and so $\Theta^{(=)}_{n,1}$ is equal to the number of prime parking functions of length $n$, which is $(n-1)^{n-1}$.
\end{proof}

An immediate corollary of the above proposition concerns complete parking preference (defined in Section \ref{prel}).

\begin{corollary}\label{complete_preference}
	The number of complete parking preferences of length $n\geq 2$ is $(n-1)^{n-1}$.
\end{corollary}

\begin{proof}
	The number of complete parking preferences of length $n$ is $\Theta^{(=)}_{n,n-1}$. Then the thesis follows from Lemma \ref{symmetry} and the above proposition.
\end{proof}
  
Referring again to Table \ref{theta=}, its partial row sums are clearly the the coefficients $\Theta^{(\leq )}_{n,k}$. Also in this case there are no correspondences in \cite{Sl}, except for $\Theta^{(\leq )}_{n,n-1}$, which seems to match sequence A071720. Such a sequence, whose closed form is $(n-1)(n+1)^{n-2}$, counts the number of spanning trees in the complete graph on $n+1$ vertices minus an edge. It would be nice to have a proof of this fact, which is still missing.

\bigskip

The results found in the present section allows us to efficiently compute several values of the coefficients $T_{n,k}$, which count permutation-invariant $k$-Naples parking functions of length $n$. They are reported in Table \ref{Tnk}, where each row is eventually constant (and equal to its rightmost entry, which is the total number of parking preferences of the relevant length).

\begin{center}
	\begin{table}[h!]
		\begin{tabular}{c|ccccccc}
			$T_{n,k}$ & $k=1$ & $k=2$ & $k=3$ & $k=4$ & $k=5$ & $k=6$ & $k=7$ \\
			\hline
			$n=2$ & 4       &          &          &          &          &          & \\
			$n=3$ & 23      & 27       &          &          &          &          & \\
			$n=4$ & 192     & 229      & 256      &          &          &          & \\
			$n=5$ & 2077    & 2558     & 2869     & 3125     &          &          & \\
			$n=6$ & 27808   & 35154    & 40000    & 43531    & 46656    &          & \\
			$n=7$ & 444411  & 572470   & 662519   & 726668   & 776887   & 823543   & \\
			$n=8$ & 8266240 & 10815697 & 12693504 & 14055341 & 15097600 & 15953673 & 16777216 \\ 	
		\end{tabular}
		\caption{The first values of the coefficients $T_{n,k}$.}\label{Tnk}
	\end{table}
\end{center}

The entries of this table also do not appear in \cite{Sl}, except of course for the diagonals $T_{n,n-1}=n^n$ and $T_{n,n-2}=n^n-(n-1)^{n-1}$, counting parking preferences and noncomplete parking preferences, respectively.

\section{Complete Naples parking functions}\label{enum_compl}

The relevance of complete parking preferences in characterizing Naples parking functions (see \cite{FV}) and in counting the permutation-invariant ones (as shown in the previous section) leads us naturally to address the problem of their enumeration. To this aim, we need to introduce a few further notations.

\bigskip

Given $n,k,m\geq 1$, denote with $v_{n,k}(m)$ the number of $k$-Naples parking functions $\alpha$ of length $n$ such that $u_{\alpha}(i)\geq 1$ for all $m+1\leq i\leq n$, and $u_{\alpha}(i)\geq 0$ for all $1\leq i\leq m$. Moreover, $v^{(0)}_{n,k}(m)$ counts the special subset of the above described preferences with the additional constraint that $u_{\alpha}(m)=0$. 

Notice that, since $u_{\alpha}(1)=0$ for all parking preferences $\alpha$, the number of complete $k$-Naples parking functions of length $n$ is given by $v_{n,k}(1)=v^{(0)}_{n,k}(1)$. Moreover, when $n\leq k+1$, Corollary \ref{complete_preference} implies that $v_{n,k}(1)=(n-1)^{n-1}$.

The next three propositions give a recursive description of the coefficients $v^{(0)}_{n,k}(1)$.

\begin{proposition}
	Given $n,k,m\geq 1$, with $n\geq k$, we have
	\begin{equation}\label{first}
		v^{(0)}_{n+1,k}(1)=\sum_{m=1}^{k}v_{n,k}(m).
	\end{equation}	
\end{proposition}

\begin{proof}
	Let $\alpha =(a_1 ,a_2 ,\ldots ,a_n ,a_{n+1}) \in PF_{n+1,k}$ be complete. This implies that $|\alpha|_1 =0$ (by Lemma \ref{elementary_properties}) and, using Lemma 4.2 of \cite{FV}, that $a_{n+1}\leq k+1$. We can then consider the parking preference $\tilde{\alpha}=(a_1 -1,a_2 -1,\ldots ,a_n -1)\in PP_n$ (obtained from $\alpha$ by removing the last entry and subtracting 1 to the remaining ones). Setting $m=a_{n+1}-1\leq k$, it is an easy computation to check that $u_{\tilde{\alpha}}(i)=u_{\alpha}(i+1)\geq 1$ when $i>m$, whereas $u_{\tilde{\alpha}}(i)=u_{\alpha}(i+1)-1\geq 0$ when $i\leq m$. As a consequence, $\tilde{\alpha}$ is one of the $k$-Naples parking functions of length $n$ counted by $v_{n,k}(m)$. Since the above construction of $\tilde{\alpha}$ from $\alpha$ is reversible, we get formula (\ref{first}).    
\end{proof}

\begin{proposition}
	Given $n,k,m\geq 1$, with $n\geq k$, we have
	\begin{equation}\label{second}
		v_{n,k}(m)=\sum_{i=1}^{m}v^{(0)}_{n,k}(i).
	\end{equation}
\end{proposition}

\begin{proof}
	Let $\alpha =(a_1 ,a_2 ,\ldots ,a_n ) \in PF_{n,k}$ be such that $u_{\alpha}(i)\geq 1$ for all $m+1\leq i\leq n$, and $u_{\alpha}(i)\geq 0$ for all $1\leq i\leq m$. Then either $u_{\alpha}(m)=0$ (and there are $v^{(0)}_{n,k}(m)$ such preferences) or ($m\geq 2$ and) $u_{\alpha}(m)\geq 0$ (and there are $v_{n,k}(m-1)$ such preferences). We thus have that
	\[
		v_{n,k}(m)=v^{(0)}_{n,k}(m)+v_{n,k}(m-1),
	\]
	and iterating the above formula gives (\ref{second}).
\end{proof}

\begin{proposition}
	Given $n,k,m\geq 1$, with $n\geq k$, we have
	\begin{equation}\label{third}
		v^{(0)}_{n,k}(m)={n\choose m-1}m^{m-2}v^{(0)}_{n-m+1,k}(1).
	\end{equation}
\end{proposition}

\begin{proof}
	Let $\alpha =(a_1 ,a_2 ,\ldots ,a_n ) \in PF_{n,k}$ be such that $u_{\alpha}(i)\geq 1$ for all $m+1\leq i\leq n$, and $u_{\alpha}(i)\geq 0$ for all $1\leq i\leq m$, with the additional requirement that $u_{\alpha}(m)=0$. Setting $J=\{i\leq n\, |\, a_i \geq m\}$, by Lemma \ref{shifts_and_restrictions} we get that $\tau_{m-1}(\alpha_{|_{J}})$ is a $k$-Naples parking function, which is complete since $u_{\tau_{m-1}(\alpha_{|_{J}})}(i)=u_{\alpha}(i+m-1)\geq 1$ for all $2\leq i\leq n-m+1$.
	
	Moreover, setting as usual $J^c =[1,n]\setminus J$, we have that $u_{\alpha_{|_{J^c}}}(i)=u_{\alpha}(i)\geq 0$ for all $1\leq i\leq m-1$. Let $\alpha_{|_{J^c}}=(\hat{a}_1 ,\hat{a}_2 ,\ldots ,\hat{a}_{m-1})$, and define $\beta=(m-\hat{a}_1 ,m-\hat{a}_2 ,\ldots ,m-\hat{a}_{m-1})$. Clearly $|\beta |_i =|\alpha_{|_{J^c}}|_{m-i}$, and so, for all $2\leq j\leq m-1$ (using Lemma \ref{elementary_properties})
	\begin{align*}
		u_{\beta}(j)&=\sum_{i=j}^{m-1}|\beta |_i -((m-1)-j+1)=\sum_{i=j}^{m-1}|\alpha_{|_{J^c}}|_{m-i}-((m-1)-j+1) \\
		&=\sum_{i=1}^{m-j}|\alpha_{|_{J^c}}|_i -((m-1)-j+1)=-u_{\alpha_{|_{J^c}}}(m-j+1)\leq 0,
	\end{align*}
	i.e. $\beta$ is a parking function of length $m-1$. Also, notice that the above combinatorial construction (from $\alpha_{|_{J^c}}$ to $\beta$) is clearly an involution.
	
	Summing up, we have thus proved that every preference $\alpha$ as above can be uniquely determined by choosing the subset $J$ (of cardinality $n-m+1$) from $[1,n]$, then choosing a complete $k$-Naples parking function $\tau_{m-1}(\alpha_{|_{J}})$ of length $n-m+1$ and a (classical) parking function $\beta$ of length $m-1$. This gives precisely formula (\ref{third}).  
\end{proof}

We are now ready to state our final recursion.

\begin{theorem}\label{enum_complete}
	Let $n\geq 2$ and $k\geq 1$. The number of complete $k$-Naples parking functions of length $n$ satisfies the following recursion:
	\begin{equation}\label{last}
		v^{(0)}_{n,k}(1)=
		\left\{   
		\begin{aligned}
			&(n-1)^{n-1} &\quad &n\leq k+1, \\
			&\sum_{i=1}^{k}{n-1\choose i-1}i^{i-2}(k-i+1)v^{(0)}_{n-i,k}(1) &\quad &n\geq k+2.
		\end{aligned}
		\right. 
	\end{equation}
\end{theorem}

\begin{proof}
	We have already observed that, when $n\leq k+1$, $v^{(0)}_{n,k}(1)=v_{n,k}(1)=(n-1)^{n-1}$. When $n\geq k+2$, formulas (\ref{first}) and (\ref{second}) give
	\[
		v^{(0)}_{n,k}(1)=\sum_{m=1}^{k}\sum_{i=1}^{m}v^{(0)}_{n-1,k}(i)=\sum_{i=1}^{k}(k-i+1)v^{(0)}_{n-1,k}(i),
	\]
	and we conclude using formula (\ref{third}). 
\end{proof} 

Also in the case of complete Naples parking functions, our recurrence allow us to efficiently find many of the coefficients $v^{(0)}_{n,k}(1)$. Some of them are recorded in Table \ref{vnk} (where, as in Table \ref{Tnk}, rows are eventually constant).

\begin{center}
	\begin{table}[h!]
		\begin{tabular}{c|cccccccc}
			$v^{(0)}_{n,k}(1)$ & $k=1$ & $k=2$ & $k=3$ & $k=4$ & $k=5$ & $k=6$ & $k=7$ & $k=8$ \\
			\hline
			$n=2$ & 1 &      &        &         &         &         &          & \\
			$n=3$ & 1 & 4    &        &         &         &         &          & \\
			$n=4$ & 1 & 11   & 27     &         &         &         &          & \\
			$n=5$ & 1 & 38   & 131    & 256     &         &         &          & \\
			$n=6$ & 1 & 131  & 783    & 1829    & 3125    &         &          & \\
			$n=7$ & 1 & 490  & 5136   & 15634   & 29849   & 46656   &          & \\
			$n=8$ & 1 & 1897 & 34623  & 148321  & 332869  & 561399  & 823543   & \\
			$n=9$ & 1 & 7714 & 251817 & 1505148 & 4102015 & 7735566 & 11994247 & 16777216 \\ 	
		\end{tabular}
		\caption{The first values of the coefficients $v^{(0)}_{n,k}(1)$.}\label{vnk}
	\end{table}
\end{center}

Sequence $v^{(0)}_{n,n-2}(1)$ seems to match sequence A101334 in \cite{Sl}. This is indeed the case, as we prove in our last proposition.

\begin{proposition}
	For $n\geq 3$, we have
	\[
		v^{(0)}_{n,n-2}(1)=(n-1)^{n-1}-n^{n-2}.
	\]
\end{proposition}

\begin{proof}
	Thanks to Theorem \ref{enum_complete}, we obtain
	\[
		v^{(0)}_{n,n-2}(1)=\sum_{i=1}^{n-2}{n-1\choose i-1}i^{i-2}(n-i-1)v^{(0)}_{n-i,n-2}(1)=\sum_{i=0}^{n-3}{n-1\choose i}(i+1)^{i-1}(n-i-2)^{n-i-1}.
	\] 
	
	We now use the well-known Abel's generalization of Newton's binomial formula, which has several equivalent formulations. One of those is the following (see for instance \cite{R}):
	\[
		(z+w+m)^m =\sum_{i=0}^{m}{m\choose i}w(w+i)^{i-1}(z+m-i)^{m-i}.
	\]
	
	Setting $m=n-1$, $z=-1$ and $w=1$ gives
	\[
		(n-1)^{n-1}=\sum_{i=0}^{n-1}{n-1\choose i}(i+1)^{i-1}(n-i-2)^{n-i-1}.
	\] 
	
	As a consequence, subtracting from the above expression for $(n-1)^{n-1}$ the terms corresponding to $i=n-1$ and $i=n-2$ gives the desired formula for $v^{(0)}_{n,n-2}(1)$.	
\end{proof}

\end{document}